\documentclass[12pt,fullpage,doublespace ]{article}
\usepackage[latin1]{inputenc}
\usepackage{amsmath}
\usepackage{amsfonts}
\usepackage{amssymb}
\author{Vehbi Emrah Paksoy}
\usepackage{latexsym}
\usepackage{amsfonts}
\usepackage{amsmath}
\usepackage{xypic}
\newtheorem{df}[subsection]{\textbf{Definition }}
\newtheorem{prop}[subsection]{\textbf{Proposition }}
\newtheorem{lem}[subsection]{\textbf{Lemma }}
\newtheorem{thm}[subsection]{\textbf{Theorem }}
\newtheorem{cor}[subsection]{\textbf{Corollary }}

\newcommand{\pr}{\textbf{Proof :}}
\newcommand{\G}{$\mathbb{G}$}
\newcommand{\B}{\begin{flushright} $\Box$ \end{flushright}}
\title{Mirror Principle for Flag Manifolds}
\date{}
\begin{document}
\maketitle
\openup 1.5 \jot

\begin{abstract}
In this paper, using mirror principle developped by Lian, Liu and Yau \cite{LLY1,LLY2,LLY3, LYS,LLLY,LLY} we obtained the A and B series for the equivariant tangent bundles over homogenous spaces using Chern polynomial. This is necessary to obtain related cohomology valued series for given arbitrary vector bundle and multiplicative characteristic class. Moreover this can be used as a valuable testing ground for the theories which associates quantum cohomologies and J functions of non-abelian quotient to abelian quotients via quantization.
\footnote{I would like to thank Bong H. Lian for his precious helps and guideance}
\end{abstract}

\section{Introduction}

It is an interesting question to obtain A series for equivarant tangent bundles and Chern Polynomials since this will be necessary to obtain A series for a general vector bundle and multiplicative characteristic class. Now assume $\mathbb{T}$ is an
algebraic torus and $X$ be a $\mathbb{T}$-manifold with a
$\mathbb{T}$ equivariant embedding in $Y:= \mathbb{P}^{m_{1}}\times
\cdots \times \mathbb{P}^{m_{l}}$ such that pull backs of hyperplane
classes $H=(H_{1},\ldots, H_{l})$ generate $H^{2}(X,\mathbb{Q})$. We
will use the same notations for equivariant classes and their
restriction to $X$. Let $\check{K}\subset H_{2}(X)$ be the set of
points in $H_{2}(X,\mathbb{Z})_{free} $ in the dual of the closure
of the K\"ahler cone of $X$. $\check{K}$ is a semi-group and defines
a partial ordering $\succeq$ on $H_{2}(X,\mathbb{Q})_{free}$.
Explicitly $r \preceq d \mbox{ iff } d-r \in \check{K}.$ If $\{
\check{H}_{j} \} $ is the dual basis for $\{ H_{i} \} $ in
$H_{2}(X),  r\preceq d \Leftrightarrow d-r=d_{1}\check{H}_{1}+\cdots
+ d_{l}\check{H}_{l}$ where $d_{i}, i=1,\ldots , l$ are nonnegative
integers. Let $X=Fl(n)$ be the complete flag variety. The first Chow ring $A^{1}(X)\cong H^{2}(X, \mathbb{Z})$ is generated by $\mathfrak{S}_{i}=c_{1}(L_{\lambda_{i}}), i=1,\ldots, n-1$ and $\lambda_{i}$ is the dominant weight of torus action with $\lambda_{i}=(1, \ldots ,1,0, \ldots ,0)$ first $i$ terms are 1's. Here $L_{\lambda_{i}}$ is the line bundle over $X$ associated the 1 dimensional representation with respect to weight $\lambda_{i}$. For more on the homogenous manifolds one can consult \cite{BR,F2,F1,MO,TU}.
\section{Basics of Mirror Principle}
We will define stable pointed map moduli for a
general projective $\mathbb{T}$ space $X$ with $\mathbb{T}$ is an
algebraic torus acting on $X$. Let  $M_{0,n}(d,X)$ be the degree
$d\in A_{1}(X)$ arithmetic genus 0, $n-$pointed stable map moduli
stack with target $X$(see \cite{FP}, \cite{K}). Following
\cite{LLY3}, we will not use the bar notation for compactification.
A typical element can be represented by $(C, f, x_{1},\ldots
,x_{n})$. This moduli space has a "virtual" fundamental class
$[M_{0.n}(d,X)]$ of dimension $dim X+ \langle c_{1}(X),d \rangle
+n-3$ similar to the fundamental class in topology. For more details
for constructions see \cite{LT}.

Let $V$ be a vector bundle on $X$. It induces a vector bundle
$V_{d},  d\in A_{1}(X)$ on $M_{0,n}(d,X)$ whose fiber at
$(C,f,x_{1},\ldots ,x_{n})$ is given by $H^{0}(C,f^{*}V)\oplus
H^{1}(C,f^{*}V)$. Another important construction is the graph space
$M_{d}(X)$ for a projective $\mathbb{T}$ manifold $X$. $M_{d}(X)$ is
the moduli stack of degree $(1,d)$ arithmetic genus 0, 0-pointed
stable maps with target $\mathbb{P}^{1}\times X$. The standard
action of $\mathbb{C}^{*}$ on $\mathbb{P}^{1}$ together with the
action of $\mathbb{T}$ on $X$ induces an action of
$\mathbb{G}=\mathbb{T}\times \mathbb{C}^{*}$  on $M_{d}(X)$. We will
denote \G -equivariant virtual class by $[M_{d}(X)] \in
A^{\mathbb{G}}_{*}(M_{d}(X))$ which has dimension $\langle
c_{1}(X),d \rangle +dim X$.

$\mathbb{C}^{*}$ fixed points of $M_{d}(X)$ plays an important
role and will be described as
\begin{displaymath}
F_{r}:= M_{0,1}(r,X)\times_{X}M_{0,1}(d-r,X).
\end{displaymath}
For any $(C_{1},f_{1},x_{1})\times (C_{2},f_{2},x_{2})\in
F_{r}$ we can obtain an element in $M_{d}(X)$ by gluing $C_{1}$ and
$C_{2}$ to $\mathbb{P}^{1}$ at 0 and $\infty$ respectively. New
curve $C$ will be mapped to $\mathbb{P}^{1}\times X$ as follows; Map
$\mathbb{P}^{1} $ identically $\mathbb{P}^{1}$ and contract
$C_{1},C_{2}$ to $0, \infty$. Map $C_{i}$ by $f_{i}$ and contract
$\mathbb{P}^{1}$ to the point $f_{1}(x_{1})=f_{2}(x_{2})$. This
defines an element $(C,f)\in M_{d}(X)$. Observe that
$F_{0}=M_{0,1}(d,X)=F_{d}$ but they will be imbedded in $M_{d}(X)$
in two different ways. For $F_{0}$ we glue the marked point to 0 and
glue the marked point to $\infty $ for $F_{d}$ in $\mathbb{P}^{1}$.
We will denote inclusion maps $i_{r} : F_{r} \hookrightarrow
M_{d}(X)$. Note that each $F_{r}$ has an evaluation map $e^{X}_{r}:
F_{r}\rightarrow X$ sending each point to the common image of the
marked points in $X$. Here are some other notations which will be
used.
\begin{itemize}
 \item{Let $L_{r}$ be the universal line bundle on
$M_{0,1}(r,X)$ which is the tangent line at the marked point.}
\item{We have natural forgetting and projection maps}
\begin{displaymath}
\rho:M_{0,1}(d,X)\rightarrow M_{0,0}(d,X),\qquad \nu: M_{d}(X)\rightarrow M_{0,0}(d,X) \mbox{ ,(see \cite{LLY1}, \cite{LLY2})}.
\end{displaymath}
with a commutative diagram
\begin{eqnarray}
\xymatrix{ F_{0}=M_{0,1}(d,X) \ar[r]^{\hspace{.5cm} i_{0}}  \ar[dr]^{\rho} & M_{d}(X) \ar[d]^{\nu} \\
 & M_{0,0}(d,X)  }
\end{eqnarray} 
\item{Let $\alpha$ be the weight of the standard $\mathbb{C}^{*}$ action on $\mathbb{P}^{1}$. Denote by $A_{*}^{\mathbb{T}}(X)(\alpha)$ the algebra obtained from the polynomial algebra $A_{*}^{\mathbb{T}}(X)[\alpha]$ by localizing with respect to all invertible elements. For an element $\beta \in A_{*}^{\mathbb{T}}(X)(\alpha)$ we let $\overline{\beta}$ be the class obtained by $\alpha \mapsto -\alpha$ in $\beta$. Introduce formal variables $\zeta=(\zeta_{1},\ldots, \zeta_{m})$ such that $\overline{\zeta_{i}}=-\zeta_{i}, \forall i$. Let $\mathbf{\mathcal{R}}=\mathbb{C}[\mathcal{T}^{*}][\alpha]$ where $\mathcal{T}^{*}$ is the dual of the Lie algebra of $\mathbb{T}$. When we consider a multiplicative class like the Chern polynomial $c_{\mathbb{T}}(x)=\sum_{i=0}^{r} c_{i}x^{i}$, We extend the ground field to $\mathbb{C}(x)$.}
\item{ For each $d$ let $\varphi :M_{d}(X) \rightarrow W_{d}$ be \G -equivariant map into smooth manifold( or orbifold ) $W_{d}$ such that $\mathbb{C}^{*}$ fixed point components in $W_{d}$ are \G- invariant submanifolds $Y_{r}$ satisfying $\varphi^{-1}(Y_{r})= F_{r}$. Construction of such maps and spaces are given in \cite{LLY1}, \cite{LLY2}. In particular, for a smooth manifold $X$ let
\begin{displaymath}
\tau : X \rightarrow \mathbb{P}^{m_{1}}\times \cdots \times \mathbb{P}^{m_{l}} :=Y
\end{displaymath}
be an equivariant projective embedding inducing an isomorphism $A^{1}(X) \simeq A^{1}(Y).$ Then we have a \G- equivariant embedding $M_{d}(X) \rightarrow M_{d}(Y)$ and we can construct \G- equivariant map $ M_{d}(Y) \rightarrow W_{d}:=N_{d_{1}}\times \cdots \times N_{d_{l}}$  with $N_{d_{i}} \simeq \mathbb{P}^{(m_{i}+1)d_{i}+m_{i}}$ which are linear sigma models for $\mathbb{P}^{m_{i}}(\mbox{ see }\cite{LLY1})$. Therefore we obtain a map
\begin{displaymath}
\varphi : M_{d}(X) \rightarrow W_{d}
\end{displaymath}
satisfying above condition. Let $\kappa_{a}$ be the equivariant
hyperplane class in $W_{d}$ which is pulled back from $N_{d_{a}}$
and denote the equivariant hyperplane class on $Y$ by $H_{a}$ also
pulled back from $\mathbb{P}^{m_{a}} $ to $Y$. Let $Y_{r} , 0
\preceq r \preceq d$ be $\mathbb{C}^{*}$ fixed point components of
$W_{d}$ which are \G -equivariantly isomorphic to $Y=Y_{0}$ and
$j_{r} :Y_{r} \hookrightarrow W_{d}$ be the inclusion map. We have
$j_{r}^{*}\kappa_{a}=H_{a}+\langle r, H_{a} \rangle \alpha$.}

Consider the commutative diagram
\begin{displaymath}
\xymatrix{ F_{r} \ar[d]^{e} \ar[r]^{i_{r}} & M_{d}(X) \ar[d]^{\varphi} & \\
Y_{r} \ar[r]^{j_{r}} & W_{d}. }
\end{displaymath}
Following proposition helps us to carry the computations to $W_{d}$
from $M_{d}(X)$ which is easier to deal with.
 \prop{(\cite{LLY3},
Lemma 3.2) Given $\omega \in A_{\mathbb{G}}^{*}(M_{d}(X))$ we have
the following equality on $Y_{r} \simeq Y$ for $0 \preceq r \preceq
d$.
\begin{displaymath}
\frac{j_{r}^{*}\varphi_{*}(\omega \cap [M_{d}(X)])}{e_{\mathbb{G}}(Y_{r}/W_{d}) }= e_{*}\Big ( \frac{i_{r}^{*} \cap [F_{r}]^{vir}}{e_{\mathbb{G}}(F_{r}/M_{d}(X))}\Big ).
\end{displaymath} } \B

For  $d=(d_{1}, \ldots , d_{l}) , r=(r_{1}\ldots , r_{l})
\preceq d$ we have
\begin{displaymath}
e_{\mathbb{G}}(Y_{r}/W_{d})=\prod_{a=1}^{l} \prod_{i=0}^{m_{a}} \prod_{\substack{k=0 \\ k \neq r_{a}}}^{d_{a}} (H_{a}-u_{a,i}-(k-r_{a})\alpha)
\end{displaymath}
where $u_{i,a}$ are $\mathbb{T}$ weights of $\mathbb{P}^{m_{a}}$.
\item \emph{Note that a class $\phi \in H^{2}_{\mathbb{T}}(X)$ has a \G-equivariant extension $\hat{\phi}\in H^{2}_{\mathbb{G}}(W_{d})$ determined by $j_{r}^{*}\hat{\phi}=\phi+\langle \phi,r \rangle \alpha$ by localization theorem. We denote by $\langle H^{2}_{\mathbb{T}}(X) \rangle $ the ring generated by $H^{2}_{\mathbb{T}}(X)$ and $R_{d}$ the ring generated by their lifts. So we have the following definition  from
\cite{LLY3}}

\df{Let  $\Gamma \in H^{*}_{\mathbb{T}}(Y)$. A list $P: P_{d} \in H_{\mathbb{G}}^{*}(W_{d}), d \succeq 0$ is an $\Gamma$- Euler data  on  if
\begin{displaymath}
 \Gamma \cdot  j_{r}^{*}P_{d}= \overline{j_{0}^{*}P_{r}}\cdot j_{0}^{*}P_{d-r}
\end{displaymath} }

\emph{An immediate observation is when we apply $\tau^{*}$ we get
\begin{displaymath}
\tau^{*} \Gamma \cdot \tau^{*}j_{r}^{*}P_{d} =\overline{\tau^{*}j_{0}^{*}P_{r}} \cdot \tau^{*} j_{0}^{*} P_{d-r}
\end{displaymath}}
\item{ \emph{There is an interesting construction for linear sigma models for a toric variety $X$.( see \cite{LLY2},\cite{LLY3}
)}}
\item{  \emph{Whenever $t=(t_{1},\ldots ,t_{l})$ formal variable we let $d\cdot t=\sum d_{i}t_{i}, \linebreak \kappa \cdot t=\sum \kappa_{a} t_{a}, H\cdot t= \sum
H_{a}t_{a}$.}}
\end{itemize}

Fix a $\mathbb{T}$ equivariant multiplicative class
$b_{\mathbb{T}}$ and an equivariant vector bundle $V=V^{+}\oplus
V^{-}$ where $V^\pm$ are convex and concave bundles on $X$. We will
assume $\Omega= \frac{b_{\mathbb{T}}(V^{+})}{b_{\mathbb{T}}(V^{-})}$
is a well defined class on $X$. For such a vector bundle we have
\begin{displaymath}
V_{d} \rightarrow M_{0,0}(d,X), \qquad \mathcal{U}_{d}\rightarrow M_{d}(X)
\end{displaymath}
where $\mathcal{U}_{d}=\nu^{*}V_{d}$. Define the linear maps
\begin{eqnarray}
i_{r}^{vir}: A_{\mathbb{G}}^{*}(M_{d}(X)) & \longrightarrow & A^{\mathbb{T}}_{*}(X)(\alpha) \nonumber \\
i_{r}^{vir}\omega &:= & (e^{X}_{r})_{*}( \frac{i^{*}_{r}\omega \cap [F_{r}]}{e_{\mathbb{G}}(F_{r}/M_{d}(X))})
\nonumber
\end{eqnarray}

For a given concavex bundle $V$ on $X$ and $b_{\mathbb{T}}$ we
put
\begin{eqnarray}
A^{V,b_{\mathbb{T}}}(t)=A(t) &:= & e^{-H\cdot t/ \alpha}\sum_{d}A_{d}e^{d\cdot t} \nonumber \\
A_{d} &:= & i_{0}^{vir}\nu^{*}b_{T}(V_{d}) \nonumber
\end{eqnarray}

\noindent
where $A_{0}=\Omega$ and the sum is taken over all
$d=(d_{1},\ldots ,d_{l})\in \mathbb{Z}^{l}_{+}.$ We call $A(t)$ the
$A$ series associated to $V \mbox { and } b_{\mathbb{T}}.$ In
particular if we specialize $b_{\mathbb{T}}$ to the unit class we
have
\begin{eqnarray}
\mathbb{I}(t)=e^{-H\cdot t/\alpha} \sum_{d} \mathbb{I}_{d} e^{d\cdot t}, \qquad \mathbb{I}_{d}=i_{0}^{vir}1_{d}.
\end{eqnarray}
here $1_{d}$ is the unit class in $M_{d}(X)$.

\df{ Let $\Omega \in A_{\mathbb{T}}^{*}(X) $ be invertible. We call a power series of the form
\begin{displaymath}
B(t):=e^{-H\cdot t/\alpha} \sum_{d} B_{d} e^{d\cdot t}, \qquad B_{d}\in A^{\mathbb{T}}_{*}(X)(\alpha)
\end{displaymath}
an $\Omega$ - Euler series if
\begin{displaymath}
\sum_{0\preceq r \preceq d} \int_{X} \Omega^{-1}\cap \overline{B_{r}} \cdot B_{d-r} e^{(H+r\alpha)\cdot \zeta} \in \mathbf{\mathcal{R}}[[\zeta]]
\end{displaymath}
for all $d$. }

\prop{  $A^{V,b_{\mathbb{T}}}(t) =A(t)$ is an Euler series }

\noindent
\pr{ Cf. \cite{LLY3}, corr.3.9 } \B

\thm{\label{ed}}{ (\cite{LLY3}, Thm 3.11) Let $P: P_{d}$ be an $\Gamma $Euler data. Then
\begin{displaymath}
B(t)=e^{-H\cdot t/ \alpha}\sum_{d}\tau^{*}j_{0}^{*}P_{d} \cap \mathbb{I}_{d} e^{d \cdot t}
\end{displaymath}
is an  $\tau^{*}\Gamma$ Euler series. } \B

\emph{Recall that we have a commutative diagram of maps
which read $\nu \circ i_{0}=\rho$. So we can write
\begin{displaymath}
A_{d}=(e_{0})^{X}_{*}( \frac{\rho^{*} b_{\mathbb{T}}(V_{d}) \cap [M_{0,1}(d,X)]}{e_{\mathbb{G}}(F_{0}/M_{d}(X))})
\end{displaymath}
We can also compute $e_{\mathbb{G}}(F_{r}/M_{d}(X))$ explicitely for
$0 \preceq r \preceq d$. Although the $\mathbb{G}$ equivariant Euler
class of the normal bundle of $F_{0}$ in $M_{d}(X)$, that is
$N_{F_{0}/M_{d}(X)}$, will be used mostly, following lemma gives
such a class for every $\mathbb{C}^{*}$ -fixed point component  in
$M_{d}(X)$.}
 \lem { (\cite{LLY1},\cite{LLY3})  For $r\neq 0,d$
\begin{displaymath}
e_{\mathbb{G}}(F_{r}/M_{d}(X))=\alpha(\alpha+p_{0}^{*}c_{1}(L_{r}))\alpha(\alpha -p_{\infty}^{*}c_{1}(L_{d-r}))
\end{displaymath}
For r=0,d
\begin{displaymath}
e_{\mathbb{G}}(F_{0}/M_{d}(X))=\alpha(\alpha-c_{1}(L_{d})),\qquad e_{\mathbb{G}}(F_{d}/M_{d}(X))=\alpha(\alpha+c_{1}(L_{d}))
\end{displaymath}
where $p_{0}: F_{r} \rightarrow M_{0,1}(r,X)$ and $ p_{\infty}: F_{r}\rightarrow M_{0,1}(d-r,X)$ are projections.} \B

\cor{ If we denote the degree of $\alpha$ in a class $\omega \in A_{*}^{\mathbb{T}}(X)(\alpha)$ by $deg_{\alpha}\omega$ then $deg_{\alpha}A_{d} \leq -2$ }

\noindent
\pr{ We have
\begin{displaymath}
A_{d}=(e_{0})^{X}_{*}( \frac{\rho^{*} b_{\mathbb{T}}(V_{d}) \cap [M_{0,1}(d,X)]}{e_{\mathbb{G}}(F_{0}/M_{d}(X))} )= (e_{0})^{X}_{*}( \frac{\rho^{*} b_{\mathbb{T}}(V_{d}) \cap [M_{0,1}(d,X)]}
{\alpha(\alpha-c_{1}(L_{d}))})
\end{displaymath}
by previous lemma. So $deg_{\alpha}A_{d}\leq -2$ } \B

\emph{In particular when $\mathbb{I}_{d}$ is concerned we have a
better estimate for $\alpha $ degree.}

\prop{ $\forall d, \quad deg_{\alpha} \mathbb{I}_{d} \leq min(-2, -\langle c_{1}(X),d \rangle)$ }

\noindent
\pr{ If $\langle c_{1}(X), d \rangle \leq 2 $ then previous corollary gives the result. So assume $\langle c_{1}(X),d \rangle >2.$ Recall that the class $[M_{0,1}(d,X)]$ is of dimension $s=exp.dim M_{0,1}(d,X)=\langle c_{1}(X),d \rangle +dim X-2.$ Set $c=c_{1}(L_{d})$ then $c^{k}\cap [M_{0,1}(d,X)]$ is of dimension $s-k$ and so $e_{*}(c^{k}\cap [M_{0,1}(d,X)])\in A^{\mathbb{T}}_{s-k}(X)$. But this group is zero unless $s-k \leq dim X$ hence $k \geq s-dim X=\langle c_{1}(X),d \rangle-2$. By the lemma , we have
\begin{displaymath}
\mathbb{I}_{d}=\sum_{k \geq \langle c_{1}(X), d \rangle -2} \frac{1}{\alpha^{k+2}}e_{*}(c^{k}\cap[M_{0,1}(d,X)])
\end{displaymath}
hence the proposition follows. } \B

\emph{Most of the time computing $A(t)$ directly from the definition
is quite difficult. Nevertheless, provided that some conditions are
satisfied  it is possible to compute the $A$-series up to some
special operation called "Mirror Transformation" cf.
\cite{LLY1},\cite{LLY3}. The main idea of the process is to consider
another special series which we call $B$ series and if some analytic
conditions are satisfied we can get the $A$ series from this $B$
series by mirror transform. We will now, give more explanations.}

\df{ A projective $\mathbb{T}$ manifold $X$ is called a balloon manifold if the fixed point set $X^{\mathbb{T}}$ is finite and if for $p\in X^{\mathbb{T}}$ the weights of the isotropic representation $T_{p}X$ are pairwise linearly independent. The second condition is known as GKM condition.\cite{GKM} }

\emph{We will assume that the balloon manifold has the property that
if $p,q \in X^{\mathbb{T}}$ such that $i_{p}^{*}c=i_{q}^{*}c,
\forall c\in A_{\mathbb{T}}^{1}(X)$ then $p=q$. If two fixed points
$p,q$ in $X$ are connected by a $\mathbb{T}$ invariant 2-sphere we
call the sphere a balloon and denote it by $pq$. Balloon manifolds
are examined in more detail in \cite{LLY2}}

\df{ Two Euler series $A,B$ are linked if every balloon $pq$ in $X$ and every $d=\delta[pq] \succ 0$ the function $(A_{d}-B_{d})\vert_{p} \in \mathbb{C}(\mathcal{T}^{*})(\alpha)$ is regular at $\alpha=\lambda/\delta$ where $\lambda$ is the weight on the tangent line $T_{p}(pq) \subset T_{p}X$. }

\emph{Let $B(t):=e^{-H\cdot
t/\alpha}\sum_{d}\tau^{*}j_{0}^{*}P_{d}\cap \mathbb{I}_{d}e^{d\cdot
t} $ be an $\Omega=\tau^{*} \Gamma$- Euler series obtained from a
$\Gamma$-Euler data $P:P_{d}$. Following theorem is adapted from
\cite{LLY3} (thm 4.5 and corrollary 4.6).}

\thm{Suppose that at $\alpha=\lambda/\delta $ and $F=(\mathbb{P}^{1},f_{\delta},0)\in F_{0}$ we have $i_{p}^{*}\tau^{*}j_{0}^{*}P_{d}=i_{F}^{*}\rho^{*}b_{\mathbb{T}}(V_{d})$ for all $d=\delta[pq]$. Then $B(t)$ is linked to $A^{V,b_{\mathbb{T}}}(t)$ }

\emph{Now we state a theorem which relates two Euler series in the
previous setting by what we call a mirror transform. Assume
$B(t)=e^{H\cdot t/\alpha}\sum_{d} \tau^{*}j_{0}^{*}P_{d}\cap
\mathbb{I}_{d}e^{d\cdot t}$ where $\tau^{*}j_{0}^{*}P_{d}$ satisfies
the assertion of the previous theorem. In addition assume for all
$d$ we have
\begin{displaymath}
\tau^{*}j_{0}^{*}P_{d}=\Omega \alpha^{\langle c_{1}(X),d \rangle}(a+(a^{'}+a^{''}\cdot H)\alpha^{-1}+\ldots.)
\end{displaymath}
for some $a,a^{'},a^{''}\in \mathbb{C}(\mathcal{T}^{*})$ depending on $d$. Note also that $\mathbb{I}_{d}$ can be expanded as
\begin{displaymath}
\mathbb{I}_{d}=\alpha^{-\langle c_{1}(X),d \rangle}(b+(b^{'}+b^{''}\cdot H)\alpha^{-1}+\ldots.)
\end{displaymath}
for some $b,b^{'},b^{''}\in \mathbb{C}(\mathcal{T}^{*}) $ also
depending on $d$. Then}

\thm{\label{t1}}{ Suppose $A^{V,b_{\mathbb{T}} }(t),B(t)$ are as in the previous theorem and above assumptions hold. Then there exist power series $f\in\mathbf{\mathcal{R}}[[e^{t_{1}},\ldots ,e^{t_{m}}]], g=(g_{1},\ldots ,g_{m}), g_{j}\in \mathbf{\mathcal{R}}[[e^{t_{1}},\ldots ,e^{t_{m}}]]$ without constant terms such that
\begin{displaymath}
A^{V,b_{\mathbb{T}} }(t+g)=e^{f/\alpha}B(t).
\end{displaymath} }

\noindent
\pr{see \cite{LLY3} } \B

 \emph{There is an explicit method to compute
$A(t)=A^{V,b_{\mathbb{T}} }(t)$ in full generality on any balloon
manifold $X$ for arbitrary $V, b_{\mathbb{T}}$. Computations are in
terms of some $\mathbb{T}$ representations. Observe that by the
previous theorems, it is useful to understand the structure of
$i_{F}^{*}\rho^{*}b_{\mathbb{T}}(V_{d})$ and obtain an Euler series
with satisfying the conditions specified in the previous theorem so
we can compute $A$ series up to mirror transformation. Now we will
discuss some part of the method given in \cite{LLY3} to compute
$A(t)$.}

\emph{Recall that $V_{d}\rightarrow M_{0,1}(d,X)$ is a vector bundle
with fiber at $(C,f,x)$ is given by $H^{0}(C,f^{*}V)\oplus
H^{1}(C,f^{*}V)$. Then for a vector bundle $V$ on $X$ and
$F=(\mathbb{P}^{1},f_{\delta},0), d=\delta[pq]$ we have a
$\mathbb{T}$ representation
\begin{displaymath}
i_{F}^{*}\rho^{*}(V_{d})=H^{0}(C,f^{*}V)\oplus H^{1}(C,f^{*}V)
\end{displaymath}
which is the value of $b_{\mathbb{T}} $ for a trivial bundle over a
point. So the method uses the $\mathbb{T}$ representations of
related bundles on each balloon $pq \simeq \mathbb{P}^{1}$}

\emph{Let $V$ be any $\mathbb{T}$ equivariant vector bundle on $X$
and let
\begin{displaymath}
\xymatrix{ 0 \ar[r] & V_{N} \ar[r] & \cdots \ar[r] & V_{1} \ar[r] & V \ar[r] & 0}
\end{displaymath}
be an equivariant resolution. Then by Euler-Poincar\'e Principal,
\begin{displaymath}
[H^{0}(\mathbb{P}^{1},f^{*}_{\delta}V)]-[H^{1}(\mathbb{P}^{1},f^{*}_{\delta}V)]=\sum_{a} (-1)^{a+1}([H^{0}(\mathbb{P}^{1},f^{*}_{\delta}V_{a})]-[H^{1}(\mathbb{P}^{1},f^{*}_{\delta}V_{a})])
\end{displaymath}
Now, suppose each $V_{a}$ is a direct sum of $\mathbb{T}$ equivariant line bundles. Then each summand $L$ will contribute to $[H^{0}(\mathbb{P}^{1},f^{*}_{\delta}V_{a})] -[H^{1}(\mathbb{P}^{1},f^{*}_{\delta}V_{a})]$ the representations
\begin{eqnarray}
c_{1}(L)\vert_{p}-k\lambda/\delta \quad &,& k=0,\ldots , l\delta \mbox { or } \nonumber \\
c_{1}(L)\vert_{p}+k\lambda/\delta \quad &,& k=1,\ldots , -l\delta-1 \nonumber
\end{eqnarray}}

\noindent
\emph{depending on the sign of $l=\langle c_{1},[pq]
\rangle.$ For $l \geq 0 $ we get first and for $l<0$ we have the
second kind of contribution.}

\section{A- series for Fl(n)}

 \emph{Let X=$Fl(n)$ be a complete flag vaiety. $A^{1}(X)$ is generated by $\mathfrak{S}_{i}=c_{1}(L_{\lambda_{i}}), i=1,\ldots ,n-1$ and $\lambda_{i}$ is the dominant weight $\lambda_{i}=(\underbrace{1,\ldots ,1}_{i},0, \ldots ,0)$. Note that these are Schubert polynomials. Let $d=(d_{1},\ldots , d_{n-1}) $ be a class of a curve in the K\"ahler cone. Since K\"ahler cone of $X$ is generated by  $d=\sum_{i=1}^{n-1}d_{i} \check{\mathfrak{S}_{i}}$ where $\{ \check{\mathfrak{S}_{i}} \} $ forms a dual basis for $\{\mathfrak{S}_{i} \}.$ These are the Poincar\'e duals of the Schubert polynomials. Now consider
 \begin{displaymath}
 \xymatrix{ Fl(n)  \ar[r ]^{\hspace{.4cm} \tau}   & \mathbb{P}^{m_{1}}\times \cdots \times \mathbb{P}^{m_{n-1}} \ar[r]^{ j_{0}} & W_{d}:= N_{d_{1}}\times \cdots \times N_{d_{n-1}} }
 \end{displaymath}
 where $j_{0}$ is the imbedding of $\mathbb{P}^{m_{1}}\times \cdots \times \mathbb{P}^{m_{n-1}}$ as a $\mathbb{C}^{*}$ fixed point component of $W_{d}$ and for $0 \preceq r \preceq d$, all fixed point components are $\mathbb{T}$ equivariantly isomorphic to $ \mathbb{P}^{m_{1}}\times \cdots \times \mathbb{P}^{m_{n-1}}$. $\tau$ is the Plucker embedding. Here $N_{d_{i}} \simeq \mathbb{P}^{(m_{i}+1)d_{i}+m_{i}}$ and $W_{d}$ is the linear sigma model. Finally $m_{i}= {n \choose i}-1$. Let $H_{i}$ be the equivariant hyperplane classes in $\mathbb{P}^{m_{1}}\times \cdots \times \mathbb{P}^{m_{n-1}}$. Pull back of each $H_{i}$ gives the corresponding $\mathfrak{S}_{i}$. There exists \G=$\mathbb{C}^{*}\times \mathbb{T}$ -equivariant hyperplane classes $\kappa_{i}$  in $W_{d}$ with the property that  $j_{r}^{*}\kappa_{i}=H_{i}+\langle H_{i}, r\rangle \alpha$ for $0 \preceq r \preceq d.$ By the pull back of $\tau \circ j_{0}$, these $\kappa_{i}$ are taken to $\mathfrak{S}_{i}$. Again, we are using the same notation for equivalent and ordinary cohomology. We will compute the $A$ series of $Fl(n)$ for $\mathbb{T}$ equivariant tangent bundle and Chern
 polynomial.}

 \lem{\label{l1}}{ Let $[pq]$ be a class of balloon joining $p,q$.  Then
 \begin{displaymath}
 \langle y_{a}, [pq] \rangle =\int_{[pq] \simeq \mathbb{P}^{1}}y_{a}=\left \{ \begin{array}{ll}
 1 & \mbox {if } i\geq a \\
 0 & \mbox{ if } a\neq i,j \\
 -1 & \mbox{ if } j=a
 \end{array} \right.
 \end{displaymath}
 where $p=\omega, q=\omega (ij) \in S_{n}$ are permutations representing the fixed points and $(ij)$ is a transposition. }

 \noindent
 \pr{ We know that $(pq) \simeq X^{\omega(ij)}_{\omega} \simeq \mathbb{P}^{1}$, Richardson variety and $y_{a}=c_{1}(L_{\gamma_{a}}), \gamma_{a}=(\underbrace{0,\ldots , 1}_{a},0,\ldots ,0)$ is a weight of $\mathbb{T}$. Then
 \begin{displaymath}
 \langle y_{a},[pq] \rangle =\langle c_{1}(L^{*}_{\gamma_{i}}), [pq] \rangle =\langle c_{1}(\mathcal{O}(\gamma_{a,i}- \gamma_{a,j})), [pq] \rangle =\gamma_{a,i}-\gamma_{a,j}.
 \end{displaymath}
 where $\gamma_{a,i} $ means the $i$-th entry of $\gamma_{a}$.  So considering possibilities we ob tain the lemma.} \B

  \emph{Recall that in equivariant Grothendieck group we have
 \begin{eqnarray}{\label{g}}
 [TFl(n)]= \sum_{i=1}^{n}[U^{*}_{n-1}\otimes S_{\chi_{i}}]-\sum_{i=1}^{n-1}[U^{*}_{i} \otimes U_{i}]+\sum_{i=1}^{n-2} [U^{*}_{i}\otimes U_{i+1}]
 \end{eqnarray}}

 \noindent
  \emph{We know that $L_{\chi_{i}}=U_{i}/U_{i-1}, i=1,\ldots, n-1$. Of course we are using induced bundles for $\mathbb{T}$ -action  without changing the notation. Then we have
 \begin{displaymath}
 \xymatrix{ 0 \ar[r] & \ar[r] U_{1} \ar[r] & U_{2} \ar[r] & U_{2}/U_{1}\ar[r] & 0 }
 \end{displaymath}
 short exact sequence. So in Grothendieck group $[U_{2}]=[L_{\chi_{1}}]+[L_{\chi_{2}}]$.
 We can proceed for $i=1,\ldots ,n-1$ and obtain $[U_{i}]=\sum_{j=1}^{i}[L_{\chi_{j}}]$. Since the duality of vector bundles yields an involution $[V] \mapsto [V^{*}]$ in Grothendieck group. We have $[U_{i}^{*}]=\sum_{j=1}^{i} [L^{*}_{\chi_{i}}] $. So equation (\ref{g}) can be decomposed further to be
 \begin{displaymath}
 [TFl(n)]= \sum_{i=1}^{n}\sum_{a=1}^{n-1} [L^{*}_{\chi_{a}}\otimes S_{\chi_{i}}] -\sum_{i=1}^{n-1} \sum_{1\leq a,b \leq i}[L^{*}_{\chi_{a}}\otimes L_{\chi_{b}}] +\sum_{i=1}^{n-2} \sum_{\substack{1\leq a \leq i \\ 1\leq b \leq i+1}}[L^{*}_{\chi_{a}}\otimes L_{\chi_{b}}].
 \end{displaymath}}

 \emph{So we obtained a decomposition of $\mathbb{T}$ equivariant tangent bundle into line bundles in Grothendieck group. Therefore given a balloon $pq\in Fl(n)$ and $d=\delta[pq]$ together with $F=(\mathbb{P}^{1},f_{\delta},0)\in M_{0,1}(Fl(n),d)$ we have for $V=TFl(n)$ the representation $R=[H^{0}(\mathbb{P}^{1},f_{\delta}^{*}V)]-[H^{1}(\mathbb{P}^{1},f_{\delta}^{*}V)]$ is equal to
 \begin{eqnarray}{\label{g2}}
R&=&\sum_{i=1}^{n} \sum_{a=1}^{n-1}[H^{0}(\mathbb{P}^{1},f_{\delta}^{*}(L^{*}_{\chi_{a}}\otimes S_{\chi_{i}}))]-[H^{1}(\mathbb{P}^{1},f_{\delta}^{*}(L^{*}_{\chi_{a}}\otimes S_{\chi_{i}}))]  \\
 &-& \sum_{i=1}^{n-1} \sum_{1\leq a,b \leq i} [H^{0}(\mathbb{P}^{1},f_{\delta}^{*}(L^{*}_{\chi_{a}}\otimes L_{\chi_{b}}))]-[H^{1}(\mathbb{P}^{1},f_{\delta}^{*}(L^{*}_{\chi_{a}}\otimes L_{\chi_{b}}))] \nonumber \\
 &+& \sum_{i=1}^{n-2} \sum_{\substack{1\leq a \leq i \\ 1\leq b \leq i+1}}[H^{0}(\mathbb{P}^{1},f_{\delta}^{*}(L^{*}_{\chi_{a}}\otimes L_{\chi_{b}}))]-[H^{1}(\mathbb{P}^{1},f_{\delta}^{*}(L^{*}_{\chi_{a}}\otimes L_{\chi_{b}}))] \nonumber
 \end{eqnarray}}

 \noindent
 \emph{Considering (\ref{g2}) and using the method of \cite{LLY3} we can compute $i_{\rho(F)}^{*}b_{\mathbb{T}}(V_{d})$ for equivariant Chern polynomial. We will consider three
 cases.}

 \textbf{Case 1)} $[H^{0}(\mathbb{P}^{1},f_{\delta}^{*}(L^{*}_{\chi_{a}}\otimes L_{\chi_{b}}))]-[H^{1}(\mathbb{P}^{1},f_{\delta}^{*}(L^{*}_{\chi_{a}}\otimes L_{\chi_{b}}))], 1\leq a \leq i, 1\leq b \leq i+1$

 \emph{Note $c_{1}(L^{*}_{\chi_{a}}\otimes L_{\chi_{b}})\vert_{p}=(y_{a}-y_{b})\vert_{p}=u_{\omega(a)} -u_{\omega(b)}, \omega \in S_{n}$ corresponds to $p$ and $i_{p}^{*}y_{a}=u_{\omega(a)}$(\cite{TU}). We also have
 \begin{displaymath}
 l_{ab}=\langle L^{*}_{\chi_{a}}\otimes L_{\chi_{b}} ,[pq] \rangle =\langle y_{a}-y_{b},[pq] \rangle =\langle \mathcal{O}(\lambda_{s}-\lambda_{t}),[pq] \rangle=\lambda_{s}-\lambda_{t}
 \end{displaymath}
 where $\lambda=\chi_{b}-\chi_{a}$ and $pq \simeq X^{\omega(st)}_{\omega}, q=\omega(st).$ So as in lemma (\ref{l1}) we can compute $l_{ab}$. Namely assuming $a<b$ we obtain
 \begin{displaymath}
 l_{ab}=\left \{ \begin{array}{ll}
 0 & \mbox{ if } s,t\neq a,b \\
 -1 & \mbox{ if }  t=a \mbox{ or } s=b  \\
 1 & \mbox{ if } s=a, t\neq b \mbox{ or } s=b, t\neq a \\
 2 & \mbox{ if } s=a, t=b
 \end{array} \right.
 \end{displaymath}
 and for $b<a$ we have
  \begin{displaymath}
 l_{ab}=\left \{ \begin{array}{ll}
 0 & \mbox{ if } s,t\neq a,b \\
 -1 & \mbox{ if } s=b \mbox{ and } t\neq a \mbox{ or }s=a, t\neq b \\
 1 & \mbox{ if } s=a \mbox{ or } t=b \\
 -2 & \mbox{ if } s=b, t=a
 \end{array} \right.
 \end{displaymath}
 for $1\leq a \leq i, 1\leq b \leq i+1, 1\leq i \leq n-2.$ This contributes as $(x+c_{1}(L^{*}_{\chi_{a}}\otimes L_{\chi_{b}})\vert_{p}-k\lambda/\delta) $ for $l_{ab} \geq 0, k=0,\ldots, l_{ab}\delta$ and $(x+c_{1}(L^{*}_{\chi_{a}}\otimes L_{\chi_{b}})\vert_{p}+k\lambda/\delta)$ for $l_{ab}<0, k=1,\ldots , -l_{ab}\delta-1$ and eventually we get
 \begin{displaymath}
 \frac{\prod_{l_{ab} \geq 0} \prod_{k=0}^{l_{ab}\delta}(x+c_{1}(L^{*}_{\chi_{a}}\otimes L_{\chi_{b}})\vert_{p}-k\lambda/\delta)}{\prod_{l_{ab} < 0} \prod_{k=1}^{-l_{ab}\delta-1}(x+c_{1}(L^{*}_{\chi_{a}}\otimes L_{\chi_{b}})\vert_{p}+k\lambda/\delta)}
 \end{displaymath}}

 \textbf{Case 2)} $-[H^{0}(\mathbb{P}^{1},f_{\delta}^{*}(L^{*}_{\chi_{a}}\otimes L_{\chi_{b}}))]-[H^{1}(\mathbb{P}^{1},f_{\delta}^{*}(L^{*}_{\chi_{a}}\otimes L_{\chi_{b}}))], 1\leq a,b \leq n-1$

 \emph{Similarly $i^{*}_{p}(y_{a}-y_{b})=u_{\omega(a)-\omega(b)}$ and set
 \begin{displaymath}
 l_{ab}=\langle c_{1}(L^{*}_{\chi_{a}}\otimes L_{\chi_{b}}),[pq]  \rangle
 \end{displaymath}
 which can be computed as before and we obtain
 \begin{displaymath}
 \frac{\displaystyle \prod_{l_{ab} < 0} \prod_{k=1}^{-l_{ab}\delta-1}(x+c_{1}(L^{*}_{\chi_{a}}\otimes L_{\chi_{b}})\vert_{p}+k\lambda/\delta)}{\displaystyle \prod_{l_{ab} \geq  0} \prod_{k=0}^{l_{ab}\delta}(x+c_{1}(L^{*}_{\chi_{a}}\otimes L_{\chi_{b}})\vert_{p}-k\lambda/\delta)}
 \end{displaymath}
because of the negative sign in front.}

\textbf{Case 3) } $[H^{0}(\mathbb{P}^{1},f_{\delta}^{*}(L^{*}_{\chi_{a}}\otimes S_{\chi_{i}}))]-[H^{1}(\mathbb{P}^{1},f_{\delta}^{*}(L^{*}_{\chi_{a}}\otimes S_{\chi_{i}}))], 1\leq i \leq n, 1\leq a \leq n-1$.

\emph{This time we have $c_{1}(L^{*}_{\chi_{a}}\otimes
S_{\chi_{i}})\vert_{p}=u_{\omega(a)}-u_{i}$ and
\begin{displaymath}
 l_{a} =\langle c_{1}(L^{*}_{\chi_{a}}\otimes S_{\chi_{i}}),[pq] \rangle
 \end{displaymath}
  The contribution will be
\begin{displaymath}
\frac{\displaystyle \prod_{l_{a}\geq 0} \prod_{k=0}^{l_{a}\delta}(x+c_{1}(L^{*}_{\chi_{a}}\otimes S_{\chi_{i}})\vert_{p}-k\lambda/\delta)}{\displaystyle \prod_{l_{a}<0} \prod_{i=1}^{-l_{a}-1}(x+c_{1}(L^{*}_{\chi_{a}}\otimes S_{\chi_{i}})\vert_{p}+k\lambda/ \alpha)}
\end{displaymath}}

\emph{Combining all above we obtain}

\thm{\label {t}}{ Let $X=Fl(n), F=(\mathbb{P}^{1},f_{\delta},0)$ with $d=\delta[pq]$ for a balloon $pq \subset X$ where $p=\omega, q=\omega(jn)$ and $pq \simeq X^{\omega(jn)}_{\omega}.$ Then at $\alpha=\lambda/\delta$
\begin{eqnarray}
 i^{*}_{\rho(F)}b_{\mathbb{T}}(V_{d}) &=& \prod_{i=1}^{n}\prod_{a=1}^{n-1} \frac{\displaystyle \prod_{l_{a}\geq 0} \prod_{k=0}^{l_{a}\delta}(x+c_{1}(L^{*}_{\chi_{a}}\otimes S_{\chi_{i}})\vert_{p}-k\lambda/\delta)}{\displaystyle \prod_{l_{a}<0} \prod_{i=1}^{-l_{a}-1}(x+c_{1}(L^{*}_{\chi_{a}}\otimes S_{\chi_{i}})\vert_{p}+k\lambda/ \alpha)}\nonumber \\
 &\cdot& \prod_{i=1}^{n-1} \prod_{1\leq a,b \leq i} \frac{\displaystyle \prod_{l_{ab} < 0} \prod_{k=1}^{-l_{ab}\delta-1}(x+c_{1}(L^{*}_{\chi_{a}}\otimes L_{\chi_{b}})\vert_{p}+k\lambda/\delta)}{\displaystyle \prod_{l_{ab} \geq  0} \prod_{k=0}^{l_{ab}\delta}(x+c_{1}(L^{*}_{\chi_{a}}\otimes L_{\chi_{b}})\vert_{p}-k\lambda/\delta)} \nonumber \\
 &\cdot & \prod_{i=1}^{n-2} \prod_{\substack{ 1\leq a \leq i \\  1 \leq b \leq i+1 }}  \frac{\displaystyle \prod_{l_{ab} \geq 0} \prod_{k=0}^{l_{ab}\delta}(x+c_{1}(L^{*}_{\chi_{a}}\otimes L_{\chi_{b}})\vert_{p}-k\lambda/\delta)}{\displaystyle \prod_{l_{ab} < 0} \prod_{k=1}^{-l_{ab}\delta-1}(x+c_{1}(L^{*}_{\chi_{a}}\otimes L_{\chi_{b}})\vert_{p}+k\lambda/\delta)} \nonumber
 \end{eqnarray} } \B

\emph{For $d=\sum d_{i}\check{\mathfrak{S}_{i}}$, in
$A_{\mathbb{G}}^{*}(W_{d}) $ define
\begin{eqnarray}
Q_{d} &=& \prod_{i=1}^{n}\prod_{a=1}^{n-1} \frac{\displaystyle \prod_{d_{a}-d_{a-1} \geq 0}\prod_{k=0}^{d_{a}-d_{a-1}}(x+\kappa_{a}-\kappa_{a-1}-u_{i}-k\alpha)}{\displaystyle \prod_{d_{a}-d_{a-1} <0} \prod_{k=1}^{d_{a-1}-d_{a}-1}(x+\kappa_{a}-\kappa_{a-1}-u_{i}+ k\alpha)} \nonumber \\
 &\cdot& \prod_{i=1}^{n-1} \prod_{1\leq a,b \leq i} \frac{\displaystyle \prod_{d_{ab} < 0} \prod_{k=1}^{-d_{ab}-1}(x+\kappa_{ab}+k\alpha)}{\displaystyle \prod_{d_{ab} \geq  0} \prod_{k=0}^{d_{ab}}(x+\kappa_{ab}-k \alpha)} \nonumber \\
 &\cdot & \prod_{i=1}^{n-2} \prod_{\substack{ 1\leq a \leq i \\  1 \leq b \leq i+1 }}  \frac{\displaystyle \prod_{d_{ab} \geq 0} \prod_{k=0}^{d_{ab}}(x+\kappa_{ab}-k \alpha)}{\displaystyle \prod_{d_{ab}< 0} \prod_{k=1}^{-d_{ab}-1}(x+\kappa_{ab}+k\alpha)} \nonumber
 \end{eqnarray}}

 \noindent
 \emph{where
 \begin{eqnarray}
 d_{ab}&=& \langle y_{a}-y_{b},d \rangle =\langle \mathfrak{S}_{a}-\mathfrak{S}_{a-1}-\mathfrak{S}_{b}+\mathfrak{S}_{b-1},d \rangle=d_{a}-d_{a-1}-d_{b}+d_{b-1} ,\nonumber \\
  \kappa_{ab}&=&\kappa_{a}-\kappa_{a-1}-(\kappa_{b}-\kappa_{b-1}).\nonumber
 \end{eqnarray}}

 \prop{ With the notations of the previous theorem, $i_{p}^{*}\tau^{*}j_{0}^{*}Q_{d}= i^{*}_{\rho(F)}b_{\mathbb{T}}(V_{d})$ .}

 \noindent
 \pr{ We have $d=\delta[pq], \alpha=\lambda/\delta$ and note that $d_{i}=\langle \mathfrak{S}_{i},\delta[pq] \rangle$ and $\tau^{*}j_{0}^{*}\kappa_{a}=\tau^{*}H_{i}=\mathfrak{S}_{i}$. So $i_{p}\tau^{*}j_{0}^{*}Q_{d}$ will give the same expression as theorem (\ref{t}).} \B

 \prop{ $B(t)=e^{-\mathfrak{S} \cdot t/\alpha}\sum_{d}\tau^{*}j_{0}^{*}Q_{d}\cap \mathbb{I}_{d} e^{d \cdot t}$ is an $\Omega-$ Euler series.  Here $\mathfrak{S}=(\mathfrak{S}_{1}, \ldots ,\mathfrak{S}_{n-1})$.}

 \noindent
 \pr{Equivariant Chern polynomial of $Fl(n)$ is given by
 \begin{displaymath}
 \Omega=\tau^{*}\Gamma=b_{\mathbb{T}}(V)=\frac{\displaystyle \prod_{i=1}^{n} \prod_{a=1}^{n-1}(x+y_{a}-u_{i}) \prod_{i=1}^{n-2}\prod_{\substack{ 1\leq a \leq i \\ 1\leq b \leq i+1}}(x+y_{a}-y_{b})}{\displaystyle \prod_{i=1}^{n-1} \prod_{1\leq a,b \leq i}(x+y_{a}-y_{b})}
 \end{displaymath}
 where $\Gamma$ is a $\mathbb{T}$ equivariant class in $H^{*}_{\mathbb{T}}(Y)$ given by
 \begin{displaymath}
 \Gamma=\frac{\displaystyle \prod_{i=1}^{n} \prod_{a=1}^{n-1}(x+H_{a}-H_{a-1}-u_{i}) \prod_{i=1}^{n-2}\prod_{\substack{ 1\leq a \leq i \\ 1\leq b \leq i+1}}(x+H_{ab})}{\displaystyle \prod_{i=1}^{n-1} \prod_{1\leq a,b \leq i}(x+H_{ab})}
 \end{displaymath}
 where $H_{ab}=H_{a}-H_{a-1}-(H_{b}-H_{b-1}).$

 We must show that $\Gamma \cdot j_{r}^{*}Q_{d}=\overline{j_{0}^{*}Q_{r}}\cdot j_{0}^{*}Q_{d-r} , 0 \preceq r \preceq d$. We will consider several cases.  Let's fix $1\leq i \leq n$ and $1\leq a \leq n-1$. For $0 \preceq r \preceq d,$
 \begin{itemize}
 \item{ If $d_{a}-d_{a-1} \geq 0, r_{a}-r_{a-1} \geq 0$ and $ d_{a}-d_{a-1}-(r_{a}-r_{a-1})\geq 0$ then we will have a term $\prod_{k=0}^{d_{a}-d_{a-1}}(x+\kappa_{a}-\kappa_{a-1}-u_{i}-k\alpha) $ in $Q_{d}$. Isolate $(x+H_{a}-H_{a-1}-u_{i})$, a part of $\Omega$ , to compute
 \begin{eqnarray}{\label{y1}}
 (x+H_{a}-H_{a-1}-u_{i}) \cdot  j_{r}^{*}\prod_{k=0}^{d_{a}-d_{a-1}}(x+\kappa_{a}-\kappa_{a-1}-u_{i}-k\alpha) = \nonumber \\
 (x+y_{a}-u_{i}) \prod_{k=0}^{d_{a}-d_{a-1}}(x+H_{a}-H_{a-1}-u_{i}+((r_{a}-r_{a-1}-k)\alpha)
 \end{eqnarray}

 \noindent
 On the other hand we have $r_{a}-r_{a-1} \geq 0$ and $ d_{a}-d_{a-1}-(r_{a}-r_{a-1})\geq 0$.  Consider
 \begin{eqnarray}
 \overline{j_{0}^{*}\prod_{k=0}^{r_{a}-r_{a-1}}(x+\kappa_{a}-\kappa_{a-1}-u_{i}-k\alpha) }\cdot
  j_{0}^{*}\prod_{k=0}^{(d-r)_{a}}(x+\kappa_{a}-\kappa_{a-1}-u_{i}-k\alpha) \nonumber
 \end{eqnarray}

\noindent
where $(d-a)_{a}=d_{a}-d_{a-1}-(r_{a}-r_{a-1})$. This becomes
\begin{eqnarray}{\label{y3}}
\prod_{k=0}^{r_{a}-r_{a-1}}(x+H_{a}-H_{a-1}-u_{i}+k\alpha) \prod_{k=0}^{(d-a)_{a}}(x+H_{a}-H_{a-1}-u_{i}-k\alpha)
  \end{eqnarray}

  Expanding (\ref{y1}) and comparing to (\ref{y3}) we clearly see that they are equal. This is also contained as an example of Euler data in \cite{LLY1}. }
 \item{$d_{a}-d_{a-1} \geq 0, r_{a}-r_{a-1} \geq 0 $ but $(d-a)_{a}<0$. In this case we still have (\ref{y1}) but this time we must consider
 \begin{eqnarray}{\label{y4}}
 \frac{ \displaystyle \prod_{k=0}^{r_{a}-r_{a-1}}(x+H_{a}-H_{a-1}-u_{i}+k\alpha)}{\displaystyle \prod_{k=1}^{-(d-a)_{a}-1}(x+H_{a}-H_{a-1}-u_{i}+k\alpha) }
 \end{eqnarray}

 \noindent
 Recall that $d_{a}-d_{a-1}-(r_{a}-r_{a-1})<0 \Rightarrow -(d-a)_{a}-1 < r_{a}-r_{a-1}$. Moreover expanding (\ref{y4}) we see that only remaining term is
 \begin{displaymath}
 (x+H_{a}-H_{a-1}-u_{i})(x+H_{a}-H_{a-1}-u_{i}-(d-r)_{a} \alpha)\cdots (x+H_{a}-H_{a-1}-u_{i}+(r_{a}-r_{a-1})\alpha)
 \end{displaymath}
 which is equal to (\ref{y1}). }
 \item{ $d_{a}-d_{a-1} \geq 0$ and $r_{a}-r_{a-1} <0$. In this case we have $(d-a)_{a}>0$ and obtain the equality
 \begin{eqnarray}
 (x+H_{a}-H_{a-1}-u_{i})&\cdot& \prod_{k=0}^{d_{a}-d_{a-1}}(x+H_{a}-H_{a-1}-u_{i}+(r_{a}-r_{a-1}-k)\alpha)  \nonumber \\
&=& \frac{ \displaystyle \prod_{k=0}^{(d-a)_{a}}(x+H_{a}-H_{a-1}-u_{i}-k\alpha)}{\displaystyle \prod_{k=1}^{r_{a-1}-r_{a}-1}(x+H_{a}-H_{a-1}-u_{i}-k\alpha) }\nonumber
 \end{eqnarray}

 \noindent
 Which is in fact
 \begin{eqnarray}
 (x+H_{a}-H_{a-1}-u_{i}) \cdot  j_{r}^{*}\prod_{k=0}^{d_{a}-d_{a-1}}(x+\kappa_{a}-\kappa_{a-1}-u_{i}-k\alpha)&=& \nonumber \\
 \frac{ \displaystyle j_{0}^{*}\prod_{k=0}^{(d-a)_{a}}(x+\kappa_{a}-\kappa_{a-1}- u_{i}-k\alpha)}{\displaystyle\Big( \overline{j_{0}^{*}\prod_{k=1}^{r_{a-1}-r_{a}-1}(x+\kappa_{a}-\kappa_{a-1} -u_{i}+k\alpha) \Big )}} \nonumber
 \end{eqnarray} }
 \item{ $d_{a}-d_{a-1}<0, r_{a}-r_{a-1} \geq 0$. Obviously this implies $(d-r)_{a}<0$. We compare
 \begin{eqnarray}{\label{y6}}
 \frac{(x+H_{a}-H_{a-1}-u_{i})}{\displaystyle j_{r}^{*}\prod_{k=1}^{d_{a-1}-d_{a}-1}(x+\kappa_{a}-\kappa_{a-1}-u_{i}+k\alpha)} \quad \mbox { and }
 \end{eqnarray}
 \begin{eqnarray}{\label{y7}}
 \frac{\displaystyle\Big( \overline{ j_{0}^{*}\prod_{k=0}^{r_{a}-r_{a-1}}(x+\kappa_{a}-\kappa_{a-1} -u_{i}-k\alpha) \Big )}}{\displaystyle j_{0}^{*}\prod_{k=1}^{-(d-a)_{a}-1}(x+\kappa_{a}-\kappa_{a-1} -u_{i}+k\alpha) }
 \end{eqnarray}

 \noindent
 If $r_{a}-r_{a-1}=-(d-r)_{a}-1=r_{a}-r_{a-1}-d_{a}+d_{a-1}+1$ then we have \linebreak $d_{a}-d_{a-1}=-1$ and no term on (\ref{y6})  except $(x+H_{a}-H_{a-1}-u_{i})$ appears. Clearly only the same term survives on (\ref{y7}) after cancellation. Otherwise, observing $-(d-r)_{a}-1 > r_{a}-r_{a-1}$ and expanding (\ref{y7}) accordingly we obtain the equality of (\ref{y6}) and (\ref{y7}).}
 \item{$d_{a}-d_{a-1}<0, r_{a}-r_{a-1}<0, (d-r)\geq 0$. This time we will compare
 \begin{eqnarray}{\label{y8}}
 \frac{(x+H_{a}-H_{a-1}-u_{i})}{\displaystyle \prod_{k=1}^{d_{a-1}-d_{a}-1}(x+H_{a}-H_{a-1}-u_{i}+(r_{a}-r_{a-1}+k)\alpha)} \qquad \mbox{ and }
 \end{eqnarray}
 \begin{eqnarray}{\label{y9}}
 \frac{\displaystyle \prod_{k=0}^{(d-r)_{a}}(x+H_{a}-H_{a-1}-u_{i}-k\alpha) }{\displaystyle  \prod_{k=1}^{r_{a-1}-r_{a}-1}(x+H_{a}-H_{a-1}-u_{i}-k\alpha) }
 \end{eqnarray}

\noindent
Observe that
\begin{eqnarray}
r_{a-1}-r_{a}-1-(d-r)_{a}=-1-(d_{a}-d_{a-1})=\left \{ \begin{array}{ll}
> 0 &  d_{a}-d_{a-1} <-1 \\
0 &  d_{a}-d_{a-1}=-1
\end{array} \right. \nonumber
\end{eqnarray}

\noindent
If $d_{a}-d_{a-1}=-1$, (\ref{y8}) is just $(x+H_{a}-H_{a-1}-u_{i})$ and same for (\ref{y9}). Note if $d_{a}-d_{a-1} <-1$ then after cancellations on (\ref{y9}) we obtain the equality again. Finally }
\item{ $d_{a}-d_{a-1}<0, r_{a}-r_{a-1}<0, (d-r)_{a}<0$.  Then we will have the equality of
\begin{displaymath}
\frac{(x+H_{a}-H_{a-1}-u_{i})}{\displaystyle \prod_{k=1}^{d_{a-1}-d_{a}-1}(x+H_{a}-H_{a-1}-u_{i}+(r_{a}-r_{a-1}+k)\alpha)}
\end{displaymath}
and
\begin{displaymath}
\frac{1}{\displaystyle \prod_{k=1}^{-(d-r)_{a}-1}(x+H_{a}-H_{a-1}-u_{i}+k\alpha) \cdot  \prod_{k=1}^{r_{a-1}-r_{a}-1}(x+H_{a}-H_{a-1}-u_{i}-k\alpha)}
\end{displaymath}
since $(d-r)_{a}<0 \Rightarrow d_{a-1}-d_{a} > r_{a-1}-r_{a}$ and we will obtain the term $(x+H_{a}-H_{a-1}-u_{i})$ in the first expression when $k=r_{a-1}-r_{a}$.}

\end{itemize}

To summarize we obtain for $0 \preceq r \preceq d$
\begin{displaymath}
\underbrace{\prod_{i=1}^{n}\prod_{a=1}^{n-1}(x+H_{a}-H_{a-1}-u_{i})}_{\Gamma^{1}} \cdot j_{r}^{*}Q_{d}^{1}=\overline{j_{0}^{*}Q_{r}^{1}}\cdot j_{0}^{*}Q_{d-r}^{1}
\end{displaymath}
where
\begin{displaymath}
Q_{d}^{1}=\prod_{i=1}^{n} \prod_{a=1}^{n-1} \frac{\displaystyle \prod_{d_{a}-d_{a-1} \geq 0}\prod_{k=0}^{d_{a}-d_{a-1}}(x+\kappa_{a}-\kappa_{a-1}-u_{i}-k\alpha)}{\displaystyle \prod_{d_{a}-d_{a-1} <0} \prod_{k=1}^{d_{a-1}-d_{a}-1}(x+\kappa_{a}-\kappa_{a-1}-u_{i}+ k\alpha)} .
\end{displaymath}
In fact preceding argument can easily be seen to be true for the other two parts composing $Q_{d}$. Namely once we set
\begin{displaymath}
\Gamma_{2}=\frac{1}{\prod_{i=1}^{n-1} \prod_{1\leq a,b \leq i}(x+H_{ab})}, Q_{d}^{2}=\prod_{i=1}^{n-1} \prod_{1\leq a,b \leq i} \frac{\displaystyle \prod_{d_{ab} < 0} \prod_{k=1}^{-d_{ab}-1}(x+\kappa_{ab}+k\alpha)}{\displaystyle \prod_{d_{ab} \geq  0} \prod_{k=0}^{d_{ab}}(x+\kappa_{ab}-k \alpha)}
\end{displaymath}
and
\begin{displaymath}
\Gamma_{3}=\prod_{i=1}^{n-2} \prod_{\substack{ 1 \leq a \leq i \\ 1 \leq b \leq i+1}} (x+H_{ab}), Q_{d}^{3}=\prod_{i=1}^{n-2} \prod_{\substack{ 1\leq a \leq i \\  1 \leq b \leq i+1 }}  \frac{\displaystyle \prod_{d_{ab} \geq 0} \prod_{k=0}^{d_{ab}}(x+\kappa_{ab}-k \alpha)}{\displaystyle \prod_{d_{ab}< 0} \prod_{k=1}^{-d_{ab}-1}(x+\kappa_{ab}+k\alpha).}
\end{displaymath}
where $d_{ab}$ and $\kappa_{ab}$ are as before. Note $Q_{0}=\Omega, Q_{d}=Q_{d}^{1}\cdot Q_{d}^{2}\cdot Q_{d}^{3}, \Omega=\Omega_{1} \cdot \Omega_{2} \cdot \Omega_{3}$ and combining all of above we obtain
\begin{displaymath}
\Gamma \cdot j_{r}^{*}Q_{d}=\overline{j_{0}^{*}Q_{r}}\cdot j_{0}^{*}Q_{d-r} , 0 \preceq r \preceq d.
\end{displaymath}
 This shows that the list $Q:Q_{d}$ gives an $\Gamma $-Euler data and then by theorem (\ref{ed}) we obtain the desired result. } \B

\emph{Now we want to compute the $\alpha$ degree of $Q_{d}$.
Observing closely we find that after possible cancellations are done
$Q_{d}^{2}\cdot Q_{d}^{3}$ can be written as
  \begin{displaymath}
Q_{d}^{2}\cdot Q_{d}^{3}= \prod_{i=1}^{n-1} \prod_{a=1}^{i} \Big (\frac{\displaystyle \prod_{d_{ia} <0} \prod_{k=1}^{d_{ia}-1}(x+\kappa_{ia}+k \alpha)}{\displaystyle \prod_{d_{ia} \geq 0} \prod_{k=0}^{d_{ia}}(x+\kappa_{ia}-k \alpha)} \Big)
 \end{displaymath}}
 The $\alpha$ degree of this expression is less than $\sum_{i=1}^{n-1}\sum_{a=1}^{i} (-d_{ia}-1)$. In addition,
 \begin{displaymath}
 deg_{\alpha} \Big ( \prod_{i=1}^{n}\prod_{a=1}^{n-1} \frac{\displaystyle \prod_{d_{a}-d_{a-1} \geq 0} \prod_{k=0}^{d_{a}-d_{a-1}}(x+\kappa_{a}-\kappa_{a-1}-u_{i}-k\alpha)}{\displaystyle \prod_{d_{a}-d_{a-1}<0} \prod_{k=1}^{d_{a-1}-d_{a}-1}(x+\kappa_{a}-\kappa_{a-1}+k\alpha)} \Big ) \leq nd_{n-1}.
 \end{displaymath}
 So we obtain $deg_{\alpha}Q_{d} \leq nd_{n-1}-\sum_{i=1}^{n-1}\sum_{a=1}^{i}(d_{ia}+1)$. Recall that $c_{1}(X)=2(\mathfrak{S}_{1}+\cdots + \mathfrak{S}_{n-1})$ for $X=Fl(n)$. Then
 \begin{eqnarray}{\label{eq1}}
 \langle c_{1}(X),d \rangle -deg_{\alpha}Q_{d} \geq 2\sum_{i=1}^{n-1}d_{i}-nd_{n-1}+\sum_{i=1}^{n-1}\sum_{a=1}^{i}(d_{ia}+1)
 \end{eqnarray}

 \noindent
  \emph{We know $d_{ia}=d_{i}-d_{i-1}-(d_{a}-d_{a-1})$. Then
 \begin{eqnarray}
 \sum_{i=1}^{n-1}\sum_{a=1}^{i}(d_{ia}+1)&=& \sum_{i=1}^{n-1}i(d_{i}-d_{i-1})-\sum_{i=1}^{n-1}\sum_{a=1}^{i}(d_{a}-d_{a-1}) \nonumber \\
 &=& -d_{1}-\cdots -d_{n-1}+nd_{n-1}-(d_{1}+\cdots + d_{n-1}) \nonumber \\
 &=&-2(d_{1}+\cdots +d_{n-1})+nd_{n-1} \nonumber
 \end{eqnarray}}

 \noindent
 \emph{Therefore (\ref{eq1}) becomes
 \begin{displaymath}
 \langle c_{1}(X),d \rangle -deg_{\alpha}Q_{d} \geq \sum_{i=1}^{n-1} i\geq 0.
 \end{displaymath}
 As a result we conclude that $\tau^{*}j_{0}^{*}Q_{d}$ satisfies the conditions of theorem (\ref{t1}) and we
 have}

 \thm{ Let $X=Fl(n)$ and $V=TX$ be the equivariant tangent bundle. The $A$- series $A^{V,b_{\mathbb{T}}(X)}$ with equivariant Chern polynomial $b_{\mathbb{T}}$ can be computed as
 \begin{displaymath}
 A(t+g)=e^{f/\alpha}B(t)
 \end{displaymath}
 where $B(t)=e^{-y\cdot t/\alpha}\sum_{d}\tau^{*}j_{0}^{*}Q_{d}\cap \mathbb{I}_{d}$ and $f,g$ are formal power series given as in theorem (\ref{t1}) } \B

Vehbi E. Paksoy
Claremont McKenna College
1325 N. College Ave, D 323
Claremont, CA, 91711, USA
emrah.paksoy@cmc.edu

\end{document}